\documentclass[12pt]{amsart}
\usepackage{amscd,amsthm,amssymb,amsfonts,amsmath,epsfig}
\usepackage{url}
\usepackage[arrow,matrix,tips,2cell]{xy}
\UseHalfTwocells

\theoremstyle{plain}
\newtheorem{thm}{Theorem}[section]
\newtheorem{lemma}[thm]{Lemma}
\newtheorem{prop}[thm]{Proposition}

\theoremstyle{definition}
\newtheorem*{defn}{Definition}

\newtheorem*{remark}{Remark}
\theoremstyle{remark}
\newtheorem*{ack}{Acknowledgments}
\newtheorem*{notations}{Notations}

\newcommand{\Z}{\mathbb Z}    % Integers
\newcommand{\R}{\mathbb R}    % Real
\newcommand{\C}{\mathbb C}    % Complex
\renewcommand{\O}{\mathcal O} % for structure sheaf

\renewcommand{\>}{\rangle} %\> is already defined.
\newcommand{\LL}{\mathcal L}
\newcommand{\DD}{\mathcal D}
\newcommand{\XX}{\mathcal X}
\newcommand{\IM}{{\mathrm{Im}}\,\Omega}
\newcommand{\THN}{e^{2\pi i {\,^t\!N} Z}e^{\pi i {\,^t\!N}\Omega N}}

\begin{document}  

\begin{abstract} 
We propose a generalization of the classical theta function to higher 
cohomology of the polarization line bundle on a family of complex tori with 
positive index. The constructed cocycles vary horizontally with respect to 
the (projective) flat connection on this family coming from a heat 
operator. They also possess modular properties similar to the classical 
ones.
\end{abstract}

\title{Theta functions for indefinite polarizations}
\author{Ilia Zharkov} 
\address{Mathematics Department, Duke University, Science Drive, Durham, NC 27708}
\email{zharkov@math.duke.edu}
\maketitle

\section{Introduction}

The classical theta function is a beautiful entire function of $n$ complex 
variables
$$\Theta(Z,\Omega)=\sum_{N\in\Z^n}e^{2\pi i {\,^t\!N} Z}e^{\pi i {\,^t\!N}\Omega N}, \quad Z\in \C^n, $$
where $\Omega\in M_n(\C)$ is a symmetric matrix with positive 
definite $\IM$. Considered as a section of the polarization line bundle $L$ 
on an abelian variety, it behaves very nicely in families. Namely, it 
provides the horizontal section with respect to the flat connection on the 
push-forward bundle for such a family.

However, if the polarization $L$ of the torus 
$X=\C^n/(\Omega\Z^n\oplus\Z^n)$ is not positive, which means that $\IM$ is 
not positive definite, then
the theta-series diverges. The most common way around this is to change the
complex structure of the torus such that the line bundle $L$ is positive
for this new complex structure. But the choice of this new complex 
structure is
not canonical. Rather, for each complex torus one gets a family of abelian
varieties parameterized by some Grassmannian (see, e.g., \cite{BirkLange} for more details). 
An analytic point of view (which is morally very  similar to the first one) taken in \cite{Kempf} is to represent $H^k(X,L)$ by harmonic forms with values in $L$. Another approach used in number theory for studying certain modular properties of the $(n,1)$-definite theta series (cf. \cite{GZag}, \cite{Po1}, \cite{Po2}) is take the sum over the positive cone only. 

We propose to identify the space of {\it theta forms} $H^k(X,L)$ with the 
group cohomology $H^k(\Lambda,F)$, where
$\Lambda\cong \Z^{2n}$ is the fundamental group of $X$ with the appropriate 
action on $F:=H^0(\C^n,\O_{\C^n})$, the space of holomorphic 
functions on $\C^n$. This approach has the advantage that one never 
leaves the analytic category. The constructed representatives manifestly 
satisfy the heat equation of the connection.
  
The motivation to study this problem came from the holomorphic anomaly
equation  for $N=2$ supersymmetric $\sigma$-model in 2 dimensions with target a Calabi-Yau
3-fold \cite{BCOV}. After twisting, such a theory defines a topological QFT.
Integrating  over the moduli of Riemann surfaces gives rise to a holomorphic
anomaly for the  partition functions $\mathcal F_g$ in the associated string theory. 

The holomorphic anomaly for the  generating function for all genera
$$W=\text{exp}[\sum_{g=1}^\infty  \lambda^{2g-2}\mathcal F_g]$$
has the form of a heat equation. In the B-model  Witten \cite{Witten} identifies this
heat equation with the natural projective  connection on the family of
intermediate Jacobians over the moduli of  Calabi-Yau's. Intermediate Jacobian
is not an abelian variety but rather a  complex torus with a principal
polarization of type $(1,h^{2,1})$ (the  polarization comes from the pairing on
the middle cohomology). Thus, the  question of finding a horizontal solution
arises quite naturally.  

In this paper we restrict our calculation mostly to the case of principal
polarization of index $k$ (i.e. $\IM$ has $k$ negative eigenvalues). In this
simple case $\dim H^k(X,L)=1$, and other degree cohomology vanish. In the
degenerate case there are also higher cohomology
$H^{k+i}(X,L)=H^k(X,L)\otimes  H^i(Y,\O_Y)$, where $Y$ is the "degenerate" part
of $X$, and no natural connection is guaranteed to exist. The generalization to
non-principal  polarizations is straight forward by introducing theta forms
with  characteristics.

The structure of the paper is as follows. In section 2 we review the general 
theory of the projective flat connections. We describe the connection in 
details in the particular case of a family of polarized complex tori.
In section 3 we identify the space of theta forms for principal polarizations 
with the $H^k(\Lambda,H^0(\C^n,\O_{\C^n}))$ and show that it is 1 dimensional. 
Using Koszul resolution of $\Z^{2n}$ we write down explicitly a cocycle in the 
above cohomology. In section 4 we study the action of the modular group on the 
space of theta forms. Namely, the modular group action affects differentials in 
the Koszul resolution by acting on the basis of $\Z[\Lambda]$. Establishing the 
isomorphism between cohomology $H^k(\Lambda,F)$ computed from different 
resolutions reveals modular properties of particular representatives. 
Finally, we conclude with an observation that the chosen representatives solve 
the heat equation already on the level of cocycles, and generalize the 
construction to non-principal polarizations.

It may be worth pointing out that the main calculations done in the paper have a lot in common with the classical case. Computation of $H^k(\Lambda,H^0(\C^n,\O_{\C^n}))$ and finding explicit representatives are very much in the spirit of finding the invariants of $\Lambda$-action on holomorphic functions in $\C^n$. The proofs of all auxiliary lemmas for the modular group action are essentially the same as those for the classical theta functions, and thus are not given in full details. The only major new ingredient is to use Fourier transform for non-Schwartz functions for finding  coboundaries between modular related cocycles.

\begin{notations}
${\,^t\!}A$ denotes the transpose of a matrix $A$. In the exponent $\pi i$ means 
$\pi\sqrt{-1}$. $\Z[G]$ stands for the group ring of $G$. $k\<V_1,\dots,V_s\>$, 
$k=\Z,\R$ or $\C$, denotes the $k$-linear span of $V_i$, where $V_i$ can be free 
$k$-modules or abstract elements.
\end{notations}

\begin{ack}
I am very grateful to R. Donagi, K. Fukaya, M. Movshev, T. Pantev and D. 
Tamarkin for very helpful conversations. Also I would like to thank Max-Planck 
Intitut f\"ur Mathematik in Bonn, where the major part of the work has been 
done, for its hospitality and financial support. The final stage of the work was 
supported in part by the NSF grant DMS-9729992.
\end{ack}

\section{Flat projective connections}
In this section we recall the general facts about the flat projective 
connections coming from heat operators. Most of the material here is based on the notes \cite{BeilKazhdan} and  \cite{DonPantev}.

\subsection{Differential operators and projective connections}
Let $V\rightarrow S$ be a holomorphic vector bundle over a complex manifold 
$S$. Denote by $\mathcal D_S(V)$ the sheaf of holomorphic differential 
operators on $V$. $\mathcal D_S(V)$ is naturally filtered by order
$$0\subset\mathcal D_S^0(V)\subset\mathcal D_S^1(V)\subset\dots,$$ 
with the associated quotients given by the principal symbol sequences
\begin{equation}
0\longrightarrow\mathcal D_S^{k-1}(V)\longrightarrow\mathcal 
D_S^{k}(V)\overset{\sigma_k}{\longrightarrow}S^k T_S\otimes{\rm 
End}(V)\longrightarrow 0.
\end{equation}

Recall that the Atiyah class is the extension class 
\begin{equation}
0 \longrightarrow {\rm End}(V) \longrightarrow \mathcal A(V)\longrightarrow 
T_S \longrightarrow 0
\end{equation}
from the symbol sequence for $\mathcal D_S^1(V)$. That is, $\mathcal A(V)$ 
consists of all differential operators of order $\leq 1$, whose symbol lies 
in $T_S\otimes {\rm id}\subset T_S\otimes{\rm End}(V)$.

A {\it projective connection} on $V$ is a holomorphic splitting $\nabla$ of 
the sequence
\begin{equation}
\xymatrix{
0 \ar[r] & {\rm End}(V)/\O_S \ar[r] & \mathcal A(V)/\O_S \ar[r]
& T_S \ar[r] 
\llowertwocell_\nabla{\omit}
& 0.
}
\end{equation}

It lifts to an ordinary connection if this splitting lifts to a splitting 
of the Atiyah class sequence (2).

Now let $\pi:\mathcal X\rightarrow S$ be a smooth family of compact complex 
manifolds. That means that $\pi$ is a surjective, proper holomorphic map 
between complex manifolds, which is flat and whose fibers $X_s$ are 
connected complex manifolds. Let $\LL$ be a holomorphic line bundle on 
$\XX$. 

Denote by $\DD_\pi (\LL)$ the sheaf of all differential operators on $\LL$ 
acting along the fibers of $\pi$ (the centralizer of $\pi^{-1}\O_S$ in 
$\DD_\XX$). $\DD_\pi (\LL)$ has the order filtration 
$$0\subset\DD_\pi^0(\LL)\subset\DD_\pi^1(\LL)\subset\dots,$$ 
where $\DD_\pi^k(\LL):=(\sigma_k)^{-1}(S^k T_\pi)\cap D_\pi(\LL), \ k\geq 
0$ (here we use the natural inclusions $S^k T_\pi\hookrightarrow S^k 
T_\XX$).

There is also a filtration of $\DD_\XX(\LL)$ by order along the S direction
$$
0\subset\DD_\XX^{0,S}(\LL)\subset\DD_\XX^{1,S}(\LL)\subset\dots,
$$
where we define $\DD_\XX^{k,S}(\LL)$ inductively by
\begin{gather*}
\DD_\XX^{0,S}(\LL):=\DD_\pi (\LL),\\
\DD_\XX^{k,S}(\LL):=\{D\in\DD^k_\XX\ |\ [D,f]\in\DD_\XX^{k-1,S}(\LL),\ {\rm 
for}\ f\in \pi^{-1}\O_S\}.
\end{gather*}

Again, the quotients are given by the symbol sequences
\begin{equation}
0\longrightarrow\DD_\XX^{k-1,S}(\LL)\longrightarrow\DD_\XX^{k,S}(\LL) 
\overset{\sigma_{k,S}}{\longrightarrow}\pi^*S^k 
T_S\otimes\DD_\pi(\LL)\longrightarrow 0.
\end{equation}

\begin{defn}
A {\it heat operator} on $\LL$ relative to $\pi$ is a sheaf homomorphism 
$\mathcal H: T_S\rightarrow\pi_*\DD^{1,S}(\LL)$, that fits into the 
commutative diagram
\begin{equation} 
\xymatrix{
T_S \ar[r]^(.4){\mathcal H} \ar[dr]_(.4){{\rm id}\otimes 1}  
 & \pi_*\DD^{1,S}(\LL) \ar[d]^{\sigma_{1,S}} \\
 &  T_S\otimes\DD_\pi(\LL)
}
\end{equation}
\end{defn}

\subsection{Connections from heat operators}
In some situations a heat operator on $\LL$ yields a projective flat 
connection on the bundles $R^i\pi_*\LL\rightarrow S$ (see \cite{Hi}, \cite{BeilKazhdan}). 
A projective connection on $V=R^i\pi_*\LL\rightarrow S$ is a splitting of 
the exact sequence (3). Equivalently, we need a homomorphism 
$\nabla:T_S\rightarrow\DD_S^1(V)/\O_S$ which fits into the commutative 
diagram
\begin{equation} 
\xymatrix{
T_S \ar[r]^(.4){\nabla} \ar[dr]_(.4){{\rm id}\otimes 1}  
 & \DD_S^1(V)/\O_S \ar[d]^{\sigma_{1}} \\
 &  T_S\otimes{\rm End}(V)
}
\end{equation}
where the symbol map $\sigma_1$ is well defined on the quotient, as 
$\O_S\subset {\rm Ker}(\sigma_1)$.

The differential operators $\DD_\XX(\LL)$ act on the local sections of 
$\LL$. Hence, they act on the local sections of $R^i\pi_*\LL$ for all $i$. 
In particular, we get a map $\pi_*\DD_\XX(\LL)\rightarrow\DD_S(V)$ which 
respects the filtration by order along S, and then induces a map 
$$
\phi:\pi_*\DD^{1,S}(\LL)\rightarrow\DD_S^1(V).
$$ 
A projective connection $\nabla$ which {\it comes from a heat operator} is 
the composition 
$$
T_S\longrightarrow\pi_*\DD^{1,S}(\LL)/\O_S\longrightarrow\DD_S^1(V)/\O_S.
$$

Of course, the existence of such a connection on $R^i\pi_*\LL$ depends on 
the global geometry of the family $(\XX,\LL)$. Recall that the 
Kodaira-Spencer map for this family $\alpha:T_S\longrightarrow 
R^1\pi_*\DD_\pi^1(\LL)$ is the connecting map in the long exact sequence 
associated with the push-forward of the sequence 
$$
0\longrightarrow\DD_\pi^1(\LL)\longrightarrow\DD_\XX^{1}(\LL) 
\longrightarrow \pi^*T_S\longrightarrow 0.
$$
Then one has the following
\begin{lemma}(\cite{BeilKazhdan})
Let $(\XX,\LL)$ be a smooth family as before. Suppose

a) $\pi_*\DD_\pi(\LL)=\O_S$ and

b) the composition map $T_S\overset{\alpha} {\longrightarrow} 
R^1\pi_*\DD_\pi^1(\LL) \longrightarrow R^1\pi_*\DD_\pi(\LL)$
vanishes.

Then for any $i\geq0$, the bundle $R^i\pi_*\LL\rightarrow S$ possesses a 
unique projective flat connection given by a heat operator.
\end{lemma}

The connection is constructed as follows. Pushing forward the symbol 
sequence (4) for $k=1$, and using the projection formula one gets the long 
exact sequence of derived images
$$
0\longrightarrow\pi_*\DD_\pi(\LL)\longrightarrow\pi_*\DD^{1,S}(\LL) 
\longrightarrow T_S\otimes\pi_*\DD_\pi(\LL)\overset{\beta}{\longrightarrow} 
R^1\pi_*\DD_\pi(\LL)\longrightarrow\dots
$$
Using $\pi_*\DD_\pi(\LL)=\O_S$ and the fact that the connecting map 
$\beta: T_S\rightarrow R^1\pi_*\DD_\pi(\LL)$ can be identified with the 
composition from (b), one gets an isomorphism
$$
\pi_*\DD^{1,S}(\LL)/\O_S\simeq T_S .$$
The projective heat operator is the inverse of this isomorphism.

Very often (in fact, in all known cases) the composition in (b) vanishes 
already in $\DD_\pi^2(\LL)$. That is, the composition 
$$ 
T_S\overset{\alpha}{\longrightarrow}R^1\pi_*\DD_\pi^1(\LL)\longrightarrow 
R^1\pi_*\DD_\pi^2(\LL) $$
is already 0. Equivalently, the Kadaira-Spencer map $\alpha$ factors as 
$$
T_S \longrightarrow \pi_*S^2 T_\pi \overset{\delta}{\longrightarrow} 
R^1\pi_*\DD_\pi^1(\LL), $$
where $\delta$ is the connecting map in the push-forward of the symbol 
sequence
$$
0\longrightarrow\DD_\pi^{1}(\LL)\longrightarrow\DD_\pi^2(\LL)\longrightarrow S^2 
T_\pi\longrightarrow 0.$$
It will also be so in the case we are interested in.

\subsection{Projective connection on a family of complex tori}
We are going to apply the general theory described above to the special 
case of polarized complex tori.

Let $(\XX,\LL)$ be a family of complex tori with non-degenerate 
polarization. This means that $\LL$ restricted to any fiber $X_s$ gives a 
line bundle $L_s$, such that $c_1(L_s)\in H^2(X_s,\R)$ is a non-degenerate 
2-form. To show that this family satisfy the conditions of the Lemma 2.1, 
it is enough to do it for each fiber. 
 
Consider the filtration 
$$
\DD_{X_s}(L_s)\supset\dots\supset\DD_{X_s}^2(L_s)\supset\DD_{X_s}^1(L_s) 
\supset\O_{X_s}\supset 0. $$
There is a spectral sequence associated to this filtration which abuts to 
$H^i (X_s,\DD_{X_s}(L_s))$ with $E_1^{-p,n+p}=H^n(X_s,S^pT_{X_s})$.
\begin{equation*}
\xymatrix{
\dots & \dots & \dots & \dots & \dots \\
0 \ar[r] & H^0(X_s,S^2T_{X_s}) \ar[r]^{d_1} &  H^1(X_s,T_{X_s}) 
\ar[r]^{d_1} &  H^2(X_s,\O_{X_s}) \ar[r] & 0 \\
& 0 \ar[r] & H^0(X_s,T_{X_s}) \ar[r]^{d_1}  & H^1(X_s,\O_{X_s}) \ar[r] & 0 
\\
& &  0 \ar[r] & H^0(X_s,\O_{X_s}) \ar[r] & 0 \\
& & & 0 & 
}
\end{equation*}
The differential $d_1$ is the convolution with the class 
$c_1(L_s)-\frac{1}{2}c_1(T_{X_s})=c_1(L_s)$ (\cite{Hi}, \cite{BeilKazhdan}). Hence, every row 
is just the Koszul complex. So we have 
$H^0(X_s,\DD_{X_s}(L_s))=H^0(X_s,\O_{X_s})$, or $\pi_*\DD_\pi(\LL)=\O_S$.

To show (b) we consider the symbol sequence
$$
0 \longrightarrow \O_{X_s}\longrightarrow \DD^1_{X_s}(L_s) \longrightarrow 
T_{X_s}  \longrightarrow 0 $$
The associated long exact sequence 
\begin{equation*}
\begin{split}
H^0(X_s,T_{X_s}) & \overset{c_1(L_s)}{\longrightarrow} H^1(X_s,\O_{X_s}) 
\overset{0}{\longrightarrow}H^1(X_s,\DD^1_{X_s}(L_s)) 
\longrightarrow H^1(X_s,T_{X_s})\\
                 & \overset{c_1(L_s)}{\longrightarrow}  H^2(X_s,\O_{X_s}) 
\longrightarrow \dots
\end{split}
\end{equation*}
shows that 
$$H^1(X_s,\DD^1_{X_s}(L_s))={\rm Ker} 
(H^1(X_s,T_{X_s})\overset{c_1(L_s)}{\rightarrow}  H^2(X_s,\O_{X_s})).$$
Thus, we deduce from the second Koszul row in the above spectral sequence that 
$$H^1(X_s,\DD^1_{X_s}(L_s))\simeq H^0(X_s,S^2 T_{X_s}).$$ 
This shows (b). Namely, the Kodaira-Spencer map factors 
$$
T_S \longrightarrow \pi_*S^2 T_\pi \overset{\sim}{\longrightarrow} 
R^1\pi_*\DD_\pi^1(\LL).$$

In coordinates, the most natural choice of trivialization of $\LL$ is to 
pass to the universal covers of the tori. Let $Z_i$ be the coordinates on 
the fibers $\C^n$. The base space can be parameterized by symmetric 
 complex matrices $\Omega=\{\omega_{ij}\}$ with non-degenerate imaginary part. 
Then for the theta line bundle the heat operator which gives the connection can 
be written as
\begin{equation}
\frac{\partial}{\partial \omega_{ij}}- \frac{1}{4\pi i} 
\frac{\partial^2}{\partial Z_i \partial Z_j}.
\end{equation}
For positive definite $\IM$ this is the heat equation for the classical 
theta function.

\section{Indefinite theta forms}

Let $\XX\rightarrow S$ be a family of $n$-dimensional complex tori. We 
assume that the family is parameterized by $\mathfrak H(k,n-k)$, the space of 
non-degenerate $n\times n$ complex symmetric 
matrices $\Omega$ with the quadratic form given by $\IM$ of fixed signature 
$(k,n-k)$. Here $k$, the number of negative directions, is called {\it 
index} of the tori. The complex torus over the point $\Omega$ is described 
as $X_\Omega=\C^n/(\Omega\Z^n\oplus\Z^n)$. We denote the universal cover 
map by $p: \C^n\rightarrow X_\Omega$.

Let us now define the sheaf of sections of the theta line bundle 
$L_\Omega$. A section of $L_\Omega$ over an open subset $U\subset X_\Omega$ 
is a holomorphic function $f(Z)$ on $p^{-1}(U)$ such that 
\begin{equation*}
f(Z+M+\Omega N)= e^{- 2\pi i {\,^t\!}N Z -\pi i {\,^t\!}N \Omega N} f(Z),\ {\rm all} \ \ 
 M,N\in\Z^n. 
\end{equation*}  
The theta line bundles on $X_\Omega$, put together, form a holomorphic line 
bundle $\LL$ 
on the family $\XX$.

Lets us fix $\Omega$, and drop the subscript from $X_\Omega$ and 
$L_\Omega$.
We want to consider the cohomology groups $H^i(X,L)$.

\subsection{Homological algebra}
The fundamental group of $X=\C^n/(\Omega\Z^n\oplus\Z^n)$ splits into direct 
sum $\Lambda=\Lambda_1\oplus\Lambda_2\cong\Z^n\oplus\Z^n$ from the 
definition. We can define the action of $\Lambda_1\oplus\Lambda_2$ on 
$\O_{\C^n}$ by
\begin{equation}
(M,N):f(Z)\mapsto\THN f(Z+M+\Omega N).
\end{equation}
According to Grothendieck \cite{Groth}, we can view $\O_{\C^n}$ as a 
$\Lambda$-sheaf on $\C^n$, where the group $\Lambda=\Z^{2n}$ acts on 
$(\C^n,\O_{\C^n})$ by translations on $\C^n$ and by (8) on sections. Then 
$L$ can be identified with $\O_{\C^n}^\Lambda$, the invariants of the 
$\Lambda$-sheaf $p_*\O_{\C^n}$ on $X$. Since the action is free, there is a 
spectral sequence (\cite{Groth}, Ch. V) which abuts to  $H^i(X,L)$ with $E_2^{p,q}= 
H^p(\Lambda,H^q(\C^n,\O_{\C^n}))$. Since  $H^q(\C^n,\O)=0$ for $p>0$, the 
spectral sequence degenerates at $E_2$. Thus, $H^i(X,L)$ is naturally 
isomorphic to $H^i(\Lambda,F)$, where $F$ is the space of all holomorphic 
functions on $\C^n$ with the $\Lambda$-action given by (8). In particular, 
$H^0(X,L)$ is the space of entire functions on $\C^n$ invariant with 
respect to the $\Lambda$-action.

Before proceeding further, let us make a few comments on notations. 
We will identify $\Lambda=\Lambda_1\oplus\Lambda_2$ with 
$\Z^{2n}=\Z^n\oplus\Z^n$ by fixing the reference basis of $\Lambda$ from 
the definition of $X$. Thus, we will often use the same symbol for an 
actual element in $\Lambda$, the $2n$-column of its coordinates in the 
reference basis, or for the $n$-column of its coordinates in 
cases when this element happens to be in $\Lambda_1$ or $\Lambda_2$. Also, we 
will not 
distinguish between a real symmetric $(n\times n)$ matrix and the quadratic 
form associated to it in the $\Lambda_1$ or $\Lambda_2$ part of the reference 
basis.
When it may cause a confusion, we will comment on the specific meaning.

To construct a nice representative in $H^i(\Lambda,F)$ we need to choose a 
rather special basis in $\Lambda$. Given a real non-degenerate symmetric 
$n\times n$ matrix $Q$ we refer to a basis $\{N_1,\dots,N_n,M_1,\dots,M_n\}$ of 
$\Lambda=\Lambda_1\oplus\Lambda_2$ as {\it split with respect to $Q$}, or simply 
as {\it $Q$-split}, if the quadratic forms associated to $Q$ in 
$\Lambda_1\otimes\R$ and to $Q^{-1}$ in $\Lambda_2\otimes\R$ are positive 
definite when restricted to $\Gamma_+=\Z\<N_{k+1},\dots,N_n\>\subset\Lambda_1$ 
and $\Z\<M_{k+1},\dots,M_n\>\subset\Lambda_2$ respectively. 
We will denote by $\Gamma_-\subset\Lambda_1$ the  sublattice spanned by 
$\{N_1,\dots,N_k\}$.

\begin{remark}
It seems always possible to find a basis such that $Q\leq 0$ on $\Gamma_-$. 
$Q<0$ is too strong to require with the counterexample given by the standard 
hyperbolic quadratic form $\left ( \begin{smallmatrix} 0 & 1\\1 & 0 
\end{smallmatrix}\right )$ in $\R^2$.  
\end{remark}

From now on let us fix a split basis $\{N_1,\dots,N_n,M_1,\dots,M_n\}$ with 
respect to $\IM$. We will denote by $(NM)=\left (\begin{smallmatrix} N & 
0\\0 & M \end{smallmatrix}\right )$ the $2n\times 2n$ matrix whose columns 
are $(N_j,M_i)$. It will be convenient to assume also that ${\,^t\!}N=M^{-1}$, that is 
$(NM)\in Sp(2n,\Z)$. This is possible without loss of generality because in some 
real basis $\IM=\IM^{-1}=\text{diag}(-1,\dots,-1,1,\dots,1)$. The transformation 
matrix $\left ( \begin{smallmatrix} R & 0\\ 0 & R^{T(-1)} \end{smallmatrix} 
\right )$ acts on $\IM$ and $\IM^{-1}$ as on quadratic forms, and primitive 
integral sublattices form a dense subset in $Gr_\R(n-k,n)$.

To calculate cohomology of the free abelian group $\Lambda=\Z^{2n}$ with a 
basis $\{x_i\}$ we will use the Koszul resolution
\begin{equation}
\longrightarrow\Z[\Lambda]\otimes\wedge^2W 
\stackrel{d}{\longrightarrow}\Z[\Lambda]\otimes W 
\stackrel{d}{\longrightarrow}\Z[\Lambda] 
\stackrel{\varepsilon}{\longrightarrow}\Z \longrightarrow 0
\end{equation}
Here $W$ is the $2n$-dimensional vector space spanned by 
$w_1,\dots,w_{2n}$. 
The differential is given by
\begin{equation}
d(x\otimes w_{p_1}\dots w_{p_k})= \sum_{i=1}^k 
(-1)^{j+1}x(x_{p_i}-1)\otimes w_{p_1}\dots \hat{w}_{p_i}\dots w_{p_k}.
\end{equation}
It is also convenient to decompose 
$$W=U_1\oplus\dots\oplus U_n\oplus V_1\oplus\dots\oplus V_n
$$
with the basis $\{u_1,\dots,u_n,v_1,\dots,v_n\}$ corresponding to the basis 
 $\{N_j,M_i\}$ of $\Lambda$. Then, by taking values of 
${\rm Hom}_\Lambda (\Z[\Lambda]\otimes\wedge^* W, F)$ 
on the elements $1\otimes u_{j_1}\wedge\dots\wedge u_{j_r}\wedge 
v_{i_1}\wedge\dots\wedge v_{i_s}$, the calculation of $H^i(\Lambda,F)$ 
reduces to the cohomology of the $2n$-tuple complex $\mathcal C$ of length 
2 in each direction, with every term $\mathcal 
C_{i_1,\dots,i_p;j_1,\dots,j_q}$ equal to $F$. There are $2n$ differentials 
$d_i,\delta_j:F\rightarrow F,\ i,j=1,\dots,n$ which act on $F$ by $d_i 
f=(M_i-1)\cdot f$ and $\delta_j f=(N_j-1)\cdot f$. 

\subsection{Calculation of $H^i(\Lambda,F)$}
Before calculating the cohomology of this $2n$-complex we need the 
following elementary statement (the author thanks A. Borodin for the 
elegant proof).

\begin{lemma}
Let $\sum_{n=1}^\infty a_n e^{-n^2}$ be an absolutely convergent series. 
Then $\sum_{n=1}^\infty s_n e^{-n^2}$, where $s_n=\sum_{k=1}^n a_k$ are the 
partial sums, also converges absolutely.
\end{lemma}
\begin{proof}
The statement is equivalent to the analogous statement for convergence of 
the integrals of positive functions. Given the convergence of 
$\int\limits_1^\infty f(x) e^{-x^2} dx$ we want to say the same is true for $f(x)$ 
replaced by $F(x)=\int\limits_1^x f(t) dt$. By changing the order of integration 
we can make the following estimate:
\begin{equation*}
\begin{split}
\int\limits_1^\infty \int\limits_1^x f(t) dt \ e^{-x^2} \ dx 
   =  \int\limits_1^\infty f(t)\int\limits_t^\infty e^{-x^2} \ dx \ dt 
   \leq  \int\limits_1^\infty f(t) \frac{e^{-t^2}}{t} \ dt 
   \leq \int\limits_1^\infty f(t) e^{-t^2} dt,
\end{split}
\end{equation*} 
which proves the statement.
\end{proof}

Now we are ready to do the final calculations, which gives another proof of 
the index theorem for complex tori.
\begin{thm}
The cohomology of the total complex associated to the complex $\mathcal C$ 
is $H^i(\mathcal C)=0,$ for $i\ne k$, and $H^k(\mathcal C)\cong \C$.
\end{thm}
\begin{proof} 
First, let us note that each of the differentials $d_i$ is surjective. 
Moreover, $d_i$ is also surjective when restricted to the $\mathrm{Ker} 
(d_j)$ for any $j\ne i$. This is equivalent to saying that 
$H^q((\C^*)^n,\O)=0,\ q>0$. So the cohomology of the complex $\mathcal C$ 
can be computed from the cohomology of the $n$-tuple complex with 
differentials $\delta_j$ and terms given by the space of $\Z^n$-periodic 
holomorphic functions. 

Every such function has the Fourier expansion which is convenient to write 
in the form 
\begin{equation}
  f(Z)= \sum_{K\in\Z^n} a_K \Theta_K(Z,\Omega),\quad 
\text{where}\quad\Theta_K(Z,\Omega)=e^{\pi i {\,^t\!}K \Omega K} e^{2\pi i {\,^t\!}K Z}. 
\end{equation}
The form of the Fourier coefficients is motivated by the following simple 
form of the differential
$$ \delta_j f(Z) = N_j\cdot f(Z)-f(Z)= \sum_{K\in\Z^n} (a_{K-N_j}-a_N) 
\Theta_K(Z,\Omega).
$$

The next observation is that $\delta_j$, for $k+1\leq q\leq n$ are 
surjective (and remain such on $\mathrm{Ker}(\delta_i)$ for $i\ne j$). 
Explicitly, if $f(Z)$ is given by (11), then 
$$ F(Z)= -\sum_{K\in\Z^n} (\sum_{r=0}^{r_j} a_{K-rN_j}) e^{\pi i {\,^t\!}K \Omega 
K} e^{2\pi i {\,^t\!}K Z}
$$
is also convergent by Lemma 3.1 and defines an entire function, and 
$\delta_j F(Z)=f(Z)$. Here $r_j$ is the coefficient of $K$ written in the 
basis $\{N_1,\dots,N_n\}$, i.e., $K=\sum_{i=1}^n r_i N_i$.

So the calculation now reduces to the $k$-complex with the differentials 
$\delta_1,\dots,\delta_k$ and terms in the form
$$ f(Z)=\sum_{N_-\in\Gamma_-} 
a_{N_-}\sum_{K\in\Gamma_++N_-}\Theta_K(Z,\Omega)
$$

In this complex the differentials $\delta_1,\dots,\delta_k$ are injective. 
This easily follows from the fact that vanishing of $(a_{K-N_q}-a_K)$ for 
$1\leq q\leq k$ means that $a_K$ are constant into the $N_q$ direction. 
This forces summing over negative definite sublattice. Hence the Fourier 
series for $f(Z)$ diverges.

Thus we see that the cohomology of $\mathcal C$ reduces to a cocycle 
concentrated in $\mathcal C_{1,2,\dots,k;\emptyset}$, and given there by a 
function in the form
$$ f(Z)= \sum_{N_-\in\Gamma_-} 
a_{N_-}\sum_{K\in\Gamma_++N_-}\Theta_K(Z,\Omega)
$$
with finite number of $a_{N_-}$ non zero. 
Finally, since 
$$ N_q \cdot \sum_{K\in\Gamma_+} \Theta_K(Z,\Omega) 
=\sum_{K\in\Gamma_++N_q} \Theta_K(Z,\Omega),\quad \text{for}\ q=1,\dots,k,
$$
the cohomology $H^k(\mathcal C)$ can be represented a multiple of the cocycle
\begin{equation*} 
c=\begin{cases}
\sum_{K\in\Gamma_+} \Theta_K(Z,\Omega) & \text{in}\ \mathcal 
C_{1,2,\dots,k;\emptyset}, \\
0 & \text{in all other terms}
\end{cases}
\end{equation*}
\end{proof}

\section{Modular properties}
As in the case of the classical theta functions, we define the modular group to 
be a subgroup of $Sp(2n,\Z)$:
$$
\Gamma_{1,2}=\Bigl\{
\begin{pmatrix} A & B\\C & D \end{pmatrix}\in Sp(2n,\Z)\ \arrowvert\ \text 
{diag}({\,^t\!}A C)\quad \text{and} \quad \text{diag}({\,^t\!}B D)\quad \text{are even} 
\Bigr\}.
$$ 
We will see in this section that modular properties of the theta forms come 
most naturally from the fact that one can use different bases to construct 
the cocycle in $H^k(\Lambda,F)$. Briefly put, the modular group acts on the set 
of bases of $\Lambda$. Changing the basis affects 
differentials in the Koszul resolution. Using explicit isomorphism for the two 
resolutions which differ by an element in $\Gamma_{1,2}$ from each other, we 
will show that the entire functions representing cocycles are modular transforms 
of each other with respect to that element of $\Gamma_{1,2}$. The strategy is to
define the modular group action on $H^k(\Lambda,F)$ via its action on the bases 
and on functions representing the cocycles for the corresponding resolutions, 
and check the triviality of this action only for the generators of 
$\Gamma_{1,2}$.
 
\subsection{Identification of resolutions}
Let $\{x_i\}=\{N_j,M_l\}$ and $\{x'_i\}=\{N'_j,M'_l\}$ be two bases of $\Lambda$ 
related by $x_i=\sum_j S_{ji}x'_j$, where the matrix $S$ is an element in 
$Sp(2n,\Z)$. Then there is an induced chain homomorphism $s_*$ between the two 
resolutions which is homotopy equivalent to the identity:
\begin{equation}
\begin{CD}
@>{d'}>>\Z[\Lambda]\otimes\wedge^2W @>{d'}>>\Z[\Lambda]\otimes W 
@>{d'}>>\Z[\Lambda] @>>> \Z @>>>0\\
@.  @AA{s_*}A @AA{s_*}A @AA{\text{id}}A @\vert @.\\
@>d>>\Z[\Lambda]\otimes\wedge^2W @>d>>\Z[\Lambda]\otimes W @>d>>\Z[\Lambda] 
@>>> \Z @>>> 0\  \\
\end{CD}
\end{equation}
In particular it give rise to an isomorphism $s^*$ on the cohomology 
$H^i(\Lambda,F)$. 
Writing down $s_*$ can be quite messy and non-canonical (it gives the 
canonical map only in cohomology). But we will try to be as explicit as 
possible in a few special cases. The calculations can be somewhat simplified by 
the fact that our cocycles are concentrated in $\mathcal 
C_{1,2,\dots,k;\emptyset}$. This means that we are interested only in the 
$1\otimes w_{1}\dots w_{k}$ component of $s_*(1\otimes w_{p_1}\dots w_{p_k})$. 
With this in mind, we will see below that only the $k$ top rows of $S$ remain 
relevant.

The calculation involves finding the coefficients $R_{ij}$ in the 
expressions of the elements $(x_i-1)=R_{ij}(x'_j-1)$ in the group ring 
$\Z[\Lambda]$. Given those, it is easy to see that the collection of 
$\Z[\Lambda]$-homomorphisms $s_*$ given by
\begin{equation}
g_*: 1\otimes w_{p_1}\dots w_{p_k}\mapsto \sum_{t_1,\dots,t_k} 
(-1)^{\sigma(\{t_i\})}
\det (\{R_{p_i t_j}\}_{i,j=1}^k) \otimes w_{t_1}\dots w_{t_k}
\end{equation} 
is a chain map. In the calculations of $R_{ij}$ below we will write (when 
possible) the group operation in $\Lambda$ multiplicatively to distinguish it 
from the addition in the group ring.

Of course, there are many ways to write the decompositions 
$(x_i-1)=R_{ij}(x'_j-1)$. So let us, first, show the preferred one by example. 
Let, say, $n=4$, $k=2$. Pick, for instance, $x_5={x'_1}^3{x'_2}^2{x'_3}^4x'_6$. 
Then we can write $x_5-1={x'_1}^3{x'_2}^2{x'_3}^4(x'_6-1) 
+{x'_1}^3{x'_2}^2({x'_3}^3+{x'_3}^3+{x'_3}^2+{x'_3}+1)(x'_3-1) 
+{x'_1}^3(x'_2+1)(x'_2-1)+ ({x'_1}^2+x'_1+1)(x'_1-1)$. Thus, indeed, we see that 
in this procedure the coefficient at $1\otimes u_{1}\dots u_{2}$ of 
$s_*(1\otimes w_{p_1}\otimes w_{p_2})$ depends only on the entries in the top 2 
rows of $S$.

{\bf Type I}: $S=\left (\begin{smallmatrix} A & 0\\0 & 
A^{T(-1)}\end{smallmatrix}\right )$. It is clear that there is no interaction 
between $N$ and $M$ parts. Since we are interested only in $1\otimes 
u_{1}\dots u_{k}$ component, we can forget about the $M$ part all together.

{\bf Type Ia}: $A=\left (\begin{smallmatrix} I & 0\\ * & * 
\end{smallmatrix}\right )$, where $I$ is the $(k\times k)$ identity matrix. In 
this case the only contribution to the coefficient of $1\otimes u_{1}\dots 
u_{k}$  
comes from $s_*(1\otimes u_{1}\dots u_{k})$ and equals 1.

{\bf Type Ib}: $A=E_k$, where $E_k$ is the  $(n\times n)$-matrix with 
$(E_k)_{ij}=1$ if $i=j$ or $(i,j)=(k-1,k)$, and 0 otherwise. Here 
$N_j-1=N'_j-1, \ j\neq k,$ and $N_k-1=N'_k\cdot N'_{k-1}-1=(N'_k-1)+ 
N'_k\cdot(N'_{k-1}-1)$. 
Still, after taking the wedge product the only contribution to  $1\otimes 
u_{1}\dots u_{k}$  
comes from $s_*(1\otimes u_{1}\dots u_{k})$ and equals 1.

{\bf Type Ic}: 
$A=E_{k+1}$. Similar to the above we have 
\begin{gather*}
s_*(1\otimes u_{1}\dots u_{k})=1\otimes u_{1}\dots u_{k},\\
s_*(1\otimes u_{1}\dots u_{k-1}\wedge u_{k+1})=1\otimes u_{1}\dots 
u_{k}+N'_k\otimes u_{1}\dots u_{k-1}\wedge u_{k+1}.
\end{gather*}
And no other terms contribute to $1\otimes u_{1}\dots u_{k}$.

{\bf Type II}: $S=\left (\begin{smallmatrix} I & 0\\B & I 
\end{smallmatrix}\right )$. Here $N_j=N'_j\cdot b(M')$ and $M_j=M'_j$, 
where the $b(M')\in\Z[\Lambda]$ depends only on $M'$'s. This gives
\begin{gather*}
 N_j-1=(N'_j-1)+N'\cdot(b(M')-1),\\
 M_j-1=M'_j-1.
\end{gather*}
Thus, only contribution to the coefficient of $1\otimes u_{1}\dots u_{k}$  
comes from $s_*(1\otimes u_{1}\dots u_{k})$ and equals to 1.

{\bf Type III}: $S=\left (\begin{smallmatrix} 0 & I\\-I & 
0\end{smallmatrix}\right )$. Here we have $M_j=N'_j$ and 
$N_j=(M'_j)^{-1}$, which also gives simple expressions for $R$'s:
\begin{gather*}
 M_j-1=N'_j-1\\
 N_j-1=-(M'_j)^{-1}(M'_j-1).
\end{gather*}
Thus, $s_*(1\otimes v_{1}\dots v_{k})=1\otimes u_{1}\dots u_{k}$, and this is 
the only contribution.

To distinguish between the cohomologies computed from the two resolutions 
we will denote the corresponding $2n$-complexes by $(\mathcal C,d)$ and 
$(\mathcal C',d')$, respectively.

\subsection{Modular group action}
The modular group acts both on the set of bases and on the set of holomorphic 
functions on $\mathfrak H(k,n-k)\times \C^n$. The action on the bases given in 
the coordinates by 
\begin{gather*}
g=\left (\begin{smallmatrix} A & B\\C & D \end{smallmatrix}\right ) : 
\{N_j,M_i\} \mapsto\{N^g_j,M^g_i\},\\
(N^g_1,\dots,N^g_n,M^g_1,\dots,M^g_n) =\left ( \begin{matrix} D & -C\\-B & 
A \end{matrix} \right )   (N_1,\dots,N_n,M_1,\dots,M_n).
\end{gather*}
Thus, in the reference basis the transformation matrix is just $g^{T(-1)}$. 
Writing the action in the basis $\{N_j,M_i\}$ identifies the matrix $S$ from the
previous subsection with $(NM)^{-1} \left.{\,^t\!}\left( \begin{smallmatrix} A & B\\C & D 
\end{smallmatrix}\right)\right. (NM)$.

The next lemma defines the classical action of the modular group on the set of 
functions on $\mathfrak H(k,n-k)\times \C^n$.
\begin{lemma}
The following defines a $\Gamma_{1,2}$-action
 \begin{gather*}  
g=\left (\begin{smallmatrix} A & B\\C & D \end{smallmatrix}\right ) :
    f(Z,\Omega)\mapsto f^g(Z,\Omega),\ \text{with} \\
     f^g(Z,\Omega) = \zeta(g) \det(C\Omega^g+D)^{1/2} e^{\pi i {\,^t\!}Z C 
{\,^t\!}(C\Omega^g+D) Z} f({\,^t\!}(C\Omega^g+D) Z,\Omega^g), 
\end{gather*}   
where $\zeta(g)$ is an 8-th root of unity and we denoted $\Omega^g  
=({\,^t\!}D \Omega-{\,^t\!}B)(-{\,^t\!}C \Omega+{\,^t\!}A)^{-1}$, or equivalently $\Omega  
=(A\Omega^g+B)(C\Omega^g+D)^{-1}$.
\end{lemma}
\begin{proof}
 This is quite standard (see, for example, the Mumford's book \cite{Mum} with 
some minor modifications). Let us, for instance, show that the exponential 
factor in the 
transformation of $f(Z,\Omega)$ behaves well with respect to the group 
multiplication. Let $h=\left ( \begin{smallmatrix} A' & B'\\C' & D' 
\end{smallmatrix} \right )$, then $hg= \left ( \begin{smallmatrix} A'A+B'C & 
A'B+B'D\\C'A+D'C & C'B+D'D 
\end{smallmatrix} \right )$. The exponential factor in $f^{hg}(Z,\Omega)$ is 
given by 
\begin{equation*}
\begin{split}
 &\pi i {\,^t\!}Z (C'A+D'C) {\,^t\!}((C'A+D'C)\Omega^{hg}+(C'B+D'D)) Z=\\
=&\pi i {\,^t\!}Z ((C'A+D'C)\Omega^{hg}+(C'B+D'D)) {\,^t\!}(C'A+D'C) Z.
\end{split}
\end{equation*}
The exponential factor in $(f^g)^h$ consists of two ingredients:
\begin{equation*}
\begin{split}
h:&\ \pi i {\,^t\!}Z C(C\Omega^g+D)Z \mapsto \pi i {\,^t\!}Z 
(C'\Omega^h+D')C {\,^t\!}(C\Omega^{hg}+D) {\,^t\!}(C'\Omega^h+D') Z=\\
&=\pi i {\,^t\!}Z (C'(A\Omega^{hg}+B)(C\Omega^{hg}+D)^{-1}+D')(C\Omega^{hg}+D) 
{\,^t\!}C {\,^t\!}(C'\Omega^h+D')Z=\\
&=\pi i {\,^t\!}Z((C'A+D'C)\Omega^{hg}+(C'B+D'D)) {\,^t\!}C {\,^t\!}(C'\Omega^h+D') Z,
\end{split}
\end{equation*}
and
$$ \pi i {\,^t\!}Z C' {\,^t\!}(C'\Omega^h+D') Z.
$$
Combining these together and using the identity ${\,^t\!}(C\Omega^g+D)= (-\Omega 
C+A)^{-1}$ in various forms, matching of the two exponential factors reduces to 
showing the following:
$$ {\,^t\!}C+{\,^t\!}(-\Omega(C'A+D'C)+A'A+B'C) C'={\,^t\!}(C'A+D'C)(-\Omega C'+A').
$$
The last identity follows at once using $
{\,^t\!}C' A'={\,^t\!}A' C'$ and ${\,^t\!}D' A'-{\,^t\!}B' C'=I.$
\end{proof}

Next we define the modular transformed theta function. For 
$g=\left ( \begin{smallmatrix} A & B\\C & D \end{smallmatrix} \right )$ and 
$\Theta_K(Z,\Omega)=e^{\pi i {\,^t\!}K \Omega K} e^{2\pi i {\,^t\!}K Z}$, $K\in\Z^n$, 
we define $\Theta^g_K(Z,\Omega)$ according to the action of the above lemma.
Explicitly,
$$ 
\Theta^g_K(Z,\Omega)=\zeta \det(C\Omega^g+D)^{1/2} e^{\pi i 
{\,^t\!}Z C {\,^t\!}(C\Omega^g+D) Z} e^{\pi i {\,^t\!}K \Omega^g K + 2\pi i 
{\,^t\!}K {\,^t\!}(C\Omega^g+D) Z}.$$
Using direct computations in the spirit of Mumford's book \cite{Mum}, we discover the 
following important property of the functions $\Theta^g_K(Z,\Omega)$:
$$\Theta^g_K(Z+M+\Omega N, \Omega)= e^{- 2\pi i {\,^t\!}N Z -\pi i 
{\,^t\!}N \Omega N} \Theta^g_{K+{\,^t\!}A N+{\,^t\!}C M}(Z,\Omega).$$
Rephrasing the above identity in the language of the $\Lambda$-action gives
$$(N,M)\cdot\Theta^g_K(Z,\Omega)=\Theta^g_{K+{\,^t\!}A N+{\,^t\!}C M}(Z,\Omega).
$$
In particular,  $N^g_j$ and $M^g_i$ act as follows:
\begin{gather*} 
M^g_i\cdot\Theta^g_K(Z,\Omega)=(-CM_i,AM_i)\cdot\Theta^g_{K}(Z,\Omega) 
=\Theta^g_{K+(-{\,^t\!}A C+{\,^t\!}C A)M_i}(Z,\Omega)=\Theta^g_{K}(Z,\Omega)\\
 N^g_j\cdot\Theta^g_K(Z,\Omega)=(DN_j,-BN_j)\cdot\Theta^g_{K}(Z,\Omega) 
=\Theta^g_{K+({\,^t\!}A D-{\,^t\!}C B)N_j}(Z,\Omega)=\Theta^g_{K+N_j}(Z,\Omega).
\end{gather*}

Repeating the construction of the cocycle from the previous section, we arrive 
at the 
following important observation:
\begin{prop}
If $\{N_1,\dots,N_n,M_1,\dots,M_n\}$ is split with respect to $\IM^g$, then
\begin{equation*} 
c^g=\begin{cases}
\sum_{K\in\Gamma_+} \Theta^g_K(Z,\Omega) & \text{in}\ \mathcal 
C^g_{1,2,\dots,k;\emptyset}, \\
0 & \text{in all other terms}
\end{cases}
\end{equation*}
 is a cocycle in $\mathcal C^g$. 
\end{prop}

Note that if $\{N_1,\dots,N_n,M_1,\dots,M_n\}$ (allowing permutations) is not a 
split basis with respect to $\IM^g$, then the cohomology $H^k(\mathcal C^g)$ 
cannot be represented by a cocycle concentrated in just one place. To avoid 
dealing with such situations we will consider the cocycles in $\mathcal C^{[g]}$ 
for each left coset $[g]\in\Gamma_{1,2}/(Sl(n,\Z)\cap\Gamma_{1,2})$. The 
subgroup $Sl(n,\Z)\cap\Gamma_{1,2}=\{\left ( \begin{smallmatrix} A & 0\\0 & 
A^{T(-1)} \end{smallmatrix} \right )\} \subset\Gamma_{1,2}$ can be thought of as 
the stabilizer of the fixed Lagrangian splitting 
$\Lambda=\Lambda_1\oplus\Lambda_2$. Hence, the left cosets can be identified 
with the set of such splittings. Given a splitting 
$\Lambda=\Lambda^{[g]}_1\oplus\Lambda^{[g]}_2$ we choose $g\in [g]$, such that 
the basis $\{N^g_1,\dots,N^g_n,M^g_1,\dots,M^g_n\}$ respects this splitting and 
$\{N_1,\dots,N_n,M_1,\dots,M_n\}$ is split with respect to $\IM^g$. Thus for any 
$[g]$ and a choice of a good representative $g$ we get a cohomology class in 
$H^k(\Lambda,F)$ which is represented by the cocycle $c^g$ in $\mathcal C^g$. 

The next statement shows that this class is, in fact, independent of the choice 
of $g\in [g]$.
\begin{prop} Let $g\in[\mathrm{id}]$ be such that $\{N_j,M_i\}$ is split with 
respect to $\IM^g$. Equivalently, the basis $\{N^g_j,M^g_i\}$ respects 
$\Lambda=\Lambda_1\oplus\Lambda_2$ and is $\IM$-split. Then, $c$ and $g^*(c^g)$ 
are homologous in $\mathcal C$.
\end{prop}
\begin{proof}
It is easy to see from linear algebra that the passage from $\{N_j,M_i\}$ to 
$\{N^g_j,M^g_i\}$ can be factored into a sequence of transformations of types 
(Ia), (Ib) and (Ic) such that on each intermediate step the basis remains 
$\IM$-split. More precisely, any positive primitive sublattice $\Gamma_+^g$ can 
be transformed to $\Gamma_+$ in at most $(n-k)$ steps by applying the 
transformations of type (Ic) (with possible conjugation by types (Ia) and (Ib)). 
This, in turn, is a consequence of the following assertion: given a positive 
primitive sublattice $\Gamma_+$ and a positive vector $r$, there is a positive 
primitive sublattice $\Gamma'_+$ which contains $r$ and $\dim 
(\Gamma'_+\cap\Gamma'_+)=n-k-1$. Since the primitive integral sublattices form a 
dense subset in the real Grassmannian of positive $(n-k)$-subspaces, it is 
enough to show the last statement over $\R$. Then we can take 
$\Gamma'_+=\<r,\Gamma_+^r\>$, the subspace spanned by $r$ and the  
$(n-k-1)$-dimensional part of $\Gamma_+$ orthogonal to $r$. Obviously, the 
quadratic form restricted to $\Gamma'_+$ is positive definite. 

According to our previous calculations, the transformations of types (Ia) and 
(Ib) have no visible impact neither on the cocycle $c$ nor on the complex 
$\mathcal C$. Hence, we are reduced to show the statement of the proposition for 
the transformation of type (Ic).

Recall that for $g$ of the type (Ic) we have $N^g_j=N_j$ for all $j$, except for 
$N^g_{k+1}=N_{k+1}-N_k$. Let $\Gamma_+^g\subset\Lambda_1$ denote the positive 
sublattice of the basis $\{N^g_1,\dots,N^g_n\}$. Next we consider the following 
entire function
\begin{equation}
f=\sum_{r\geq 0}\sum_{K\in\Gamma_+^g+rN_k}\Theta_K- \sum_{r\geq 
0}\sum_{K\in\Gamma_++rN_k}\Theta_K.
\end{equation}
Though each of the two sums are divergent, their difference is the convergent 
sum over the lattice points in between $\Gamma_+$ and $\Gamma^g_+$. On the 
picture below, which shows the $\<N_k,N_{k+1}\>$ slice of the lattice, those are
represented by bold dots. The hollow dots represent terms taken with negative 
sign and $\IM$ is negative in the shaded area.

\begin{figure}[htb]
\begin{center}
  \epsfig{file=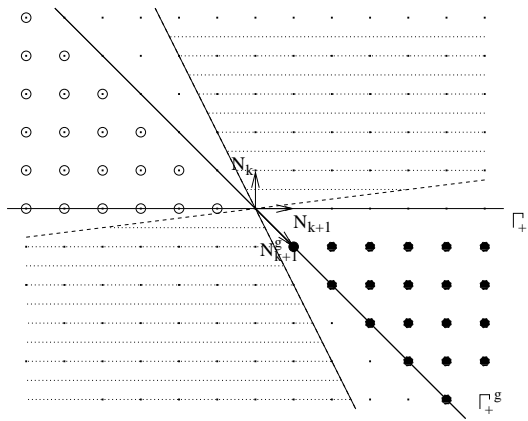}
\end{center}
\caption{The lattice points in  $\Z\<N_k,N_{k+1}\>$ for the function $f$.}
%\label{fig2}
\end{figure}

The cochain represented by this function in $\mathcal 
C_{1,2,\dots,k-1;\emptyset}$ is a coboundary between $c$ and $g^*(c^g)$. Indeed,
all differential vanish when applied to $f$ except for $(N_k-1)$ and 
$(N_{k+1}-1)$: 
\begin{gather*}
(N_k-1)\cdot f=\sum_{K\in\Gamma_+}\Theta_K-\sum_{K\in\Gamma_+^g}\Theta_K \quad 
\text{in}\ \mathcal C_{1,2,\dots,k;\emptyset},\\
(N_{k+1}-1)\cdot f=-\sum_{K\in\Gamma_+^g}\Theta_K\quad \text{in}\ \mathcal 
C_{1,2,\dots,k-1,k+1;\emptyset}.
\end{gather*}
From the calculation of (Ic) we recall that $g^*(c^g)$ has components only in 
$\mathcal C_{1,2,\dots,k;\emptyset}$ and $\mathcal 
C_{1,2,\dots,k-1,k+1;\emptyset}$ both given by $\sum_{K\in\Gamma_+} 
\Theta^g_K=\sum_{K\in\Gamma_+^g} \Theta_K$. This completes the proof.
\end{proof}
\begin{remark}
The function $f$ can be considered as a secondary object. In fact, it has a life
of its own. Similar functions appeared, for instance, in the work of G\"ottsche 
and Zagier \cite{GZag} in 
connection with the wall crossing phenomenon for Donaldson invariants. In 
general, the theta sum over the points in some simplicial $k$-cone can be 
thought of as an $k$-ary object. These "multi" theta functions appear in 
Fukaya's study of $m_k$ products among the affine Lagrangian submanifolds in a 
symplectic $2n$-torus \cite{Fuk}. It seems interesting to explore this connection 
further.
\end{remark}

The above Proposition and Lemma 4.1 allow us to define the action of the modular 
group on $H^k(\Lambda,F)$ via its action on $H^k(\mathcal C^{[g]})$. An element 
$g\in\Gamma_{1,2}$ acts on $H^k(\Lambda,F)$ as follows. Pick an element 
$h\in[g]$ such that $\{N^{}_j,M^{}_i\}$ is $\IM^h$-split. Then define
\begin{gather*}
h:\ H^k(\mathcal C)\rightarrow H^k(\mathcal C^{h})\\
h\cdot [c] \mapsto [c^{h}],
\end{gather*}
where we canonically identify $H^k(\mathcal C)$ and $H^k(\mathcal C^{h})$ with 
$H^k(\Lambda,F)$.

Finally, we are ready to prove the main result of this section.
\begin{thm}
The above action is trivial on $H^k(\Lambda,F)$. 
\end{thm}
\begin{proof}
Since $\Gamma_{1,2}$ defines an action, we need to check the triviality for the 
set of generators only.

Case 1: $g=\left ( \begin{smallmatrix} {\,^t\!}A & 0\\0 & A^{-1} \end{smallmatrix} 
\right )$. The action is defined exactly modulo the transformations of this 
type. Hence in this case the statement holds trivially.

Case 2: $g=\left (\begin{smallmatrix} I & B\\0 & I \end{smallmatrix}\right )$. 
Here we have $\Omega^g=\Omega-B$, hence $\IM=\IM^g$. Also the functions 
representing the cocycles $c$ and $c^g$ are identical: $\sum_{K\in\Gamma_+} 
\Theta_K(Z,\Omega)=\sum_{K\in\Gamma_+} \Theta^g_K(Z,\Omega)$. The matrix for the 
resolution identification is 
$$
S=(NM)\left (\begin{smallmatrix} I & 0\\{\,^t\!}B & I \end{smallmatrix}\right ) 
(NM)^{-1}= \left (\begin{smallmatrix} I & 0\\M {\,^t\!}B N^{-1} & I 
\end{smallmatrix}\right ). $$
Using the calculation for the transformation of type II we conclude that  
$g^*(c^g)=c$ already on the level of cocycles.

Case 3: $g=\left (\begin{smallmatrix} 0 & -I\\I & 0\end{smallmatrix}\right 
)$.
This seems to be the most interesting case. As for the classical theta 
functions the proof involves application of Fourier transform and Poisson 
summation formula. We will demonstrate the basic idea on the simple example 
of an 1-dimensional torus $X=\C/\tau\Z+\Z$, with $\mathrm{Im}\,\tau<0$ ($X$ 
is not an elliptic curve since the polarization is negative). According to 
the matching map on resolutions of type III we have to show that two cocycles of 
the double complex $\mathcal C$ 
\begin{equation}
\begin{CD}
0 @>{\mathbf{1}-1}>>0\\
@V{\boldsymbol{\tau}-1}VV @VVV\\
e^{\pi i \tau k^2+2\pi i k z} @>>>0 \\
\end{CD}
\qquad \text{and}\qquad
\begin{CD}
0 @>{\mathbf{1}-1}>>\zeta\frac{1}{\sqrt{\tau}}e^{-\pi i (z+n)^2/\tau}\\
@VV{\boldsymbol{\tau}-1}V @VVV\\
0 @>>>0 \\
\end{CD}
\end{equation}  
are homologous in $H^1(\mathcal C)$. To see this we define the entire 
function $f\in\mathcal C_{0;0}$ by
\begin{equation}
f(z,\tau)=\int\limits_{\mathrm{Re}\,y=k-\frac{1}{2}}\frac{e^{\pi i \tau 
y^2}e^{2\pi i (z+n)y}} {e^{2\pi i y}-1}dy.
\end{equation}
Then
\begin{equation}
\begin{split}
(\mathbf{1} -1)\cdot f(z,\tau) &=f(z+1,\tau)-f(z,\tau)=\\
&=\int\limits_{\mathrm{Re}\,y=k-\frac{1}{2}}\frac{e^{\pi i \tau y^2}e^{2\pi i 
(z+n)y}(e^{2\pi i y}-1)} {e^{2\pi i y}-1}\,dy= 
\zeta\frac{1}{\sqrt{\tau}}e^{-\pi i (z+n)^2/\tau}
\end{split}
\end{equation}
and
\begin{equation}
\begin{split}
(\boldsymbol{\tau} -1)\cdot f(z,\tau) 
&=e^{\pi i \tau+2\pi i z} f(z+\tau,\tau)-f(z,\tau)=\\
&= \int\limits_{\mathrm{Re}\,y=k-\frac{1}{2}}\frac{e^{\pi i \tau (y+1)^2}e^{2\pi i   
(z+n)(y+1)}} {e^{2\pi i y}-1}\,dy - \int\limits_{\mathrm{Re}\,y=k-\frac{1}{2}} 
\frac{e^{\pi i \tau y^2}e^{2\pi i (z+n)y}} {e^{2\pi i y}-1}\,dy \\
&= -2\pi i\: \mathrm{Res}_{y=k} (\frac{e^{\pi i \tau y^2}e^{2\pi i (z+n)y}} 
{e^{2\pi i y}-1})
=-e^{\pi i \tau k^2+2\pi i k z}.
\end{split}
\end{equation}
The minus sign appearing here reflects the choice of orientation on the line 
$\mathrm{Re}\,y=k$. This ambiguity is suppressed by the choice of 8-th root of 
unity for the modular transform. With all that, we see that the difference of 
the above cocycles is the coboundary of the function $f(z,\tau)$.

To deal with the higher dimensional case, first, we need to make some comments 
on the choice of the original split basis $\{N_j,M_i\}$. Namely, using the 
transformation of the Case 1 above, we may assume without loss of generality 
that the reference basis is $\IM$-split, that is $\{N_j,M_i\}$ can be chosen to 
be the reference basis.

In such a basis the matrix $S$ has the form $S={\,^t\!}g=\left (\begin{smallmatrix} 0 
& I\\-I & 0\end{smallmatrix}\right )$. So we can apply calculation of type III 
to conclude that $g^*(c^g)$ is represented by $\sum_{K\in\Gamma_+}\Theta_K^g$ in 
$\mathcal C_{\emptyset; 1,2,\dots,k}$.

Next, we need to choose a totally real vector subspace in 
$\C^n=\C\<N_1,\dots,N_n\>$ such that the imaginary part of $\Omega$ restricted 
to $V$ is positive definite. By induction we may also assume that $V$ has a 
codimension one filtration $V=V_0\supset V_1\supset V_2\supset\dots\supset 
V_k=\Gamma_+\otimes\R$, such that $V_i\otimes\C=\C\<N_{i+1},\dots,N_n\>$, 
$i=1,\dots,k$. The convenient property of this filtration is that $V_i$ can be 
cut out by $i$ equations $\{{\,^t\!}N_j Y=0, j=1,\dots,i\}$.

Following the analogy with the example above we define the Fourier transform of 
the function $e^{\pi i {\,^t\!}Y \Omega Y}$ by using the partial Vick rotation:
$$
\hat{f}(X)=\int\limits_{V} f(Y)e^{2\pi i {\,^t\!}Y X} dY.
$$
Then one can use the usual theory of Fourier transform for positive definite 
$\IM$. In particular, the operators of translation and multiplication by the 
corresponding character interchange. Also, using the real coordinates on 
$V=A\cdot\R^n$ one has
$$
\int\limits_V e^{\pi i {\,^t\!}Y \Omega Y} dY =\int\limits_{\R^n} e^{\pi i {\,^t\!}X({\,^t\!}A \Omega A)X}\, (\det 
A) dX =\zeta(\det A)(\det({\,^t\!}A \Omega A))^{-\frac{1}{2}} =\zeta 
(\det\Omega)^{-\frac{1}{2}}. $$
The same identity holds if the integration cycle $V$ is shifted by any vector in 
$\C^n$.

We define a cochain $f$ in $\mathcal C^{k-1}$ by the functions $f_i\in\mathcal 
C^{k-1}_{1,2,\dots,i-1;i+1,\dots,k}$, $1\leq i\leq k$:
$$
f_i(Z)=\sum_{K\in\Gamma_+}\int\limits_{V_{i-1}-\frac{1}{2}\sum_{j=i}^k N_j} 
\frac{e^{\pi i {\,^t\!}Y\Omega Y+2\pi i {\,^t\!}Y(Z+K)}}{e^{2\pi i {\,^t\!}Y N_i}-1}\, d^{n-i+1}Y.
$$
We want to show that the total differential applied to this cochain gives the 
difference between $g^*(c^g)$ and $c$ in $\mathcal C$: 
\begin{equation*} 
\xymatrix{
 & & f_1 \ar[r]^(.15){M_1-1} \ar[d]^{N_1-1} & \zeta (\det\Omega)^{-\frac{1}{2}} 
\sum_{K\in\Gamma_+} e^{-\pi i {\,^t\!}(Z+K)\Omega^{-1} (Z+K)}\\
& \dots \ar[d]^{N_{k-1}} & \dots &\\
f_k \ar[r]^{M_k-1} \ar[d]^{N_k-1} & 0 & &\\
-\sum_{K\in\Gamma_+} e^{\pi i {\,^t\!}K\Omega K+2\pi i {\,^t\!}K Z} & & &
}
\end{equation*}

Since  $f_i\in\mathcal C^{k-1}_{1,2,\dots,i-1;i+1,\dots,k}$, the differentials 
$(M_s-1)\cdot f_i$ vanish trivially for $s=i+1,\dots,k$. Also, $f_i(Z)$ is 
clearly periodic with respect to $M_s$, $s=k+1,\dots,n$. For the rest we have:
\begin{equation*}
\begin{split}
(M_s-1)\cdot f_i 
&=f_i(Z+M_s)-f_i(Z)=\\
&=\sum_{K\in\Gamma_+} \int\limits_{V_{i-1}-\frac{1}{2}\sum_{j=i}^k N_j} \frac{e^{\pi i 
{\,^t\!}Y\Omega Y+2\pi i {\,^t\!}Y (Z+K)}(e^{2\pi i {\,^t\!}Y M_s}-1)}{e^{2\pi i {\,^t\!}Y N_i}-1}\, 
d^{n-i+1}Y.
\end{split}
\end{equation*} 
But $Y\in V_{i-1}+\frac{1}{2}\sum_{j=i}^k N_j$ means that ${\,^t\!}Y M_s=0$ for $s<i$. 
Thus the only non-trivial $d$-differential is given by
\begin{equation*}
(M_i-1)\cdot f_i= \sum_{K\in\Gamma_+}\int\limits_{V_{i-1}-\frac{1}{2}\sum_{j=i}^k N_j}  
e^{\pi i {\,^t\!}Y\Omega Y+2\pi i {\,^t\!}Y (Z+K)}\, d^{n-i+1}Y.
\end{equation*}
In particular for $i=1$ we have:
\begin{equation*}
\begin{split}
(M_1-1)\cdot f_1 
&=\sum_{K\in\Gamma_+} e^{-\pi i {\,^t\!}(Z+K)\Omega^{-1} (Z+K)} 
\int\limits_{V-\frac{1}{2}\sum_{j=1}^k N_j} e^{\pi i {\,^t\!}(Y+\Omega^{-1}(Z+K))\Omega 
(Y+\Omega^{-1}(Z+K))}\, d^n Y=\\
&=\sum_{K\in\Gamma_+} e^{-\pi i {\,^t\!}(Z+K)\Omega^{-1} (Z+K)}(\zeta
(\det\Omega)^{-\frac{1}{2}})=\sum_{K\in\Gamma_+} \Theta^g_K\quad \text{in} \ 
\mathcal C_{\emptyset;1,2,\dots,k}.\\
\end{split}
\end{equation*}

Similarly, the differentials $(N_s-1)\cdot f_i$ vanish for trivial reasons for 
$s=1,\dots,i-1$. On the other hand
\begin{equation*}
\begin{split}
N_s\cdot f_i 
& = e^{2\pi i {\,^t\!}N_s Z+\pi i {\,^t\!}N_s\Omega N_s} f_s(Z+\Omega N_s)=\\ 
&= \sum_{K\in\Gamma_+}\int\limits_{V_{i-1}-\frac{1}{2}\sum_{j=i}^k N_j} \frac{e^{\pi i 
{\,^t\!}(Y+N_s)\Omega(Y+N_s) +2\pi i {\,^t\!}(Y+N_s)(Z+K)}} {e^{2\pi i {\,^t\!}Y N_i}-1}\, 
d^{n-i+1}Y \\
&=\sum_{K\in\Gamma_+}\int\limits_{V_{i-1}+N_s-\frac{1}{2}\sum_{j=i}^k N_j} \frac{e^{\pi 
i {\,^t\!}Y\Omega Y+2\pi i {\,^t\!}Y(Z+K)}}{e^{2\pi i {\,^t\!}Y N_i}-1}\, d^{n-i+1}Y.\\ 
\end{split}
\end{equation*}
For $s=k+1,\dots,n$ the domain of integration remains the same, hence $N_s\cdot 
f_i=f_i$. Also, for $s=i+1,\dots,k$, we have ${\,^t\!}Y N_i=-\frac{1}{2}+at \neq 0$, 
where $a$ is a fixed complex number with non-zero imaginary part, and $t\in\R$. 
So the integrant is a holomorphic function in some open neighborhood of the 
subspace $\R\<V_{i-1},N_s\>$ in $\C^n$. Hence the only non-trivial 
$\delta$-differential is given by
\begin{equation*}
\begin{split}
(N_i-1)\cdot f_i 
&=\sum_{K\in\Gamma_+}(\int\limits_{V_{i-1}+\frac{1}{2}N_i-\frac{1}{2}\sum_{j=i+1}^k 
N_j} \frac{e^{\pi i {\,^t\!}Y\Omega Y+2\pi i {\,^t\!}Y(Z+K)}}{e^{2\pi i {\,^t\!}Y N_i}-1}\, 
d^{n-i+1}Y-\int\limits_{V_{i-1}-\frac{1}{2}\sum_{j=i}^k N_j}\dots)\\
&\overset{P.R.}{=}
- {\sum_{K\in\Gamma_+}\int\limits_{V_i-\frac{1}{2}\sum_{j=i+1}^{k} N_j} e^{\pi i 
{\,^t\!}Y\Omega Y+2\pi i {\,^t\!}Y(Z+K)}}\, d^{n-i}Y.\\
\end{split}
\end{equation*}
The last equality is given by the Poincar\'e residue map evaluated on the 
corresponding cycles. Thus we see that $(N_s-1)\cdot f_s+(M_{s-1}-1)\cdot 
f_{s-1}=0$. Finally, using  Poisson summation formula for the function 
$$
F(X)=e^{\pi i {\,^t\!}X \Omega X+2\pi i {\,^t\!}X Z},\ \text{with}\ X_1=\dots= X_k=0,
$$
of $(n-k)$ real variables and the lattice $\Gamma_+\subset V_+^\R$, we deduce 
that 
\begin{equation*}
\begin{split}
(N_k-1)\cdot f_k 
& = -\sum_{K\in\Gamma_+} \int\limits_{V_k} e^{\pi i {\,^t\!}Y \Omega Y+2\pi i {\,^t\!}Y Z+2\pi i 
{\,^t\!}Y K} d^{n-k}Y=\\
&=-\sum_{K\in\Gamma_+} e^{\pi i {\,^t\!}K\Omega K+2\pi i {\,^t\!}K Z}=-\sum_{K\in\Gamma_+} 
\Theta_K\quad \text{in}\ \mathcal C_{1,2,\dots,k;\emptyset}.\\
\end{split}
\end{equation*}
This completes the proof of the theorem.
\end{proof}

\section{The heat equation}
The construction in the group cohomology $H^k(\Lambda,F)$ can be easily replaced 
by the one in \v Cech cohomology. Namely, one can use the covering of $X$ by one 
open patch overlapping itself at the boundary of the fundamental domain. Then 
the condition of a $k$-cochain being a cocycle translates exactly into the 
vanishing of the differential in the group cohomology complex. One can also use 
\'etale cohomology for the \'etale cover $\C^n\rightarrow X$. One way or another 
we can identify the connection on the bundle $R^k\pi_*\LL$ with the heat 
operator (7): 
\begin{equation*}
\mathcal H=\frac{\partial}{\partial \omega_{ij}}- \frac{1}{4\pi i} 
\frac{\partial^2}{\partial Z_i \partial Z_j},
\end{equation*}
acting on the elements of $F$, the entire functions on $\C^n$.

A priori, the connection is defined only projectively. But for the families 
parameterized by $\mathfrak H(k,n-k)$ there are no monodromy obstructions. 
Hence, the projective connection can be lifted to an ordinary one. The heat 
operator $\mathcal H$ defines this ordinary connection.

Clearly, $\mathcal H(\Theta_K)=0$, hence the family of the classes 
$[c(\Omega)]\in H^k(\Lambda,F_\Omega)$ defines a horizontal section of 
$R^k\pi_*\LL$ already on the level of cocycles. Moreover, a straight forward 
calculation shows that $\mathcal H(\Theta^g_K)=0$, for any $g\in\Gamma_{1,2}$. 
So the classes $[c^g(\Omega)]\in H^k(\Lambda,F_\Omega)$ also provide a 
covariantly constant section.

From the point of view of projective connection an one-dimensional vector bundle 
is not very interesting. So one would like to generalize the construction to 
non-principal polarizations. This can be done by introducing theta forms with 
characteristics. In fact, every step in the two previous sections goes through 
with just minor modifications.

Let $X_\Omega=\C^n/\Omega\Z^n\oplus\Delta\Z^n$ be a complex torus with 
polarization of type 
$$\Delta=\text{diag}(\delta_1,\dots,\delta_n),$$ 
with $\delta_i$ positive integers such that $\delta_1|\delta_2|\dots|\delta_n$. The theta line bundle 
$L_\Omega$ is defined as before by its sheaf of sections. Namely, a section of 
$L_\Omega$ over an open subset $U\subset X$ is a holomorphic function $f(Z)$ on 
$p^{-1}(U)$ such that 
\begin{equation*}
f(Z+\Delta M+\Omega N)= e^{- 2\pi i {\,^t\!}N Z -\pi i {\,^t\!}N\Omega N} f(Z),\ {\rm all} \ 
\  M,N\in\Z^n. 
\end{equation*} 
Again, $H^k(X_\Omega,L_\Omega)$ can be canonically identified with the group 
cohomology $H^k(\Lambda,F_\Omega)$, with the action of 
$\Lambda\simeq\Z^n\oplus\Z^n$ on $F_\Omega$ given by 
\begin{equation*}
(M,N):f(Z)\mapsto\THN f(Z+\Delta M+\Omega N).
\end{equation*} 

Now $H^k(\Lambda,F_\Omega)$ is of dimension 
$\det\Delta=\delta_1\delta_2\cdot\dots\delta_n$, and it can be represented by 
the functions in $\mathcal C^k_{1,2,\dots,k;\emptyset}$:
$$
\Theta[a](Z,\Omega)=\sum_{K\in\Gamma_+} e^{\pi i {\,^t\!}(K+a)\Omega (K+a)+2\pi i {\,^t\!}(K+a)Z}, $$
where $a$ is an element in the lattice $\Delta^{-1}\Z^n$. Note that shifting of 
$a$ by an element of $\Z^n$ changes the cocycle defined by $\Theta[a](Z,\Omega)$ at most by a coboundary. Hence, the cohomology class of $c[a]\in 
H^k(\Lambda,F_\Omega)$ is well defined for each $a\in\Delta^{-1}\Z^n/\Z^n$, and 
the collection of these form a basis. 

Clearly, $\mathcal H(\Theta[a](Z,\Omega))=0$, hence each family 
$c[a](\Omega)\in H^k(\Lambda,F_\Omega)$ defines a horizontal section of 
$R^k\pi_*\LL$ with respect to the flat connection given by the heat operator 
$\mathcal H$.

\bibliographystyle{plain}
\bibliography{theta}

\end{document}